\magnification=1200
\def\qed{\unskip\kern 6pt\penalty 500\raise -2pt\hbox
{\vrule\vbox to 10pt{\hrule width 4pt\vfill\hrule}\vrule}}
\centerline{DIFFERENTIATION OF SRB STATES FOR HYPERBOLIC FLOWS.}
\bigskip
\centerline{by David Ruelle\footnote{*}{Mathematics Dept., Rutgers University, and IHES.  91440 Bures sur Yvette, France.\break $<$ruelle@ihes.fr$>$}.}
\bigskip\bigskip\noindent
	{\leftskip=2cm\rightskip=2cm\sl {\bf Abstract.}  Let the ${\cal C}^3$ vector field ${\cal X}+aX$ on $M$ define a flow $(f^t_a)$ with an Axiom A attractor $\Lambda_a$ depending continuously on $a\in(-\epsilon,\epsilon)$.  Let $\rho_a$ be the SRB measure on $\Lambda_a$ for $(f^t_a)$.  If $A\in{\cal C}^2(M)$, then $a\mapsto\rho_a(A)$ is ${\cal C}^1$ on $(-\epsilon,\epsilon)$ and $d\rho_a(A)/da$ is the limit when $\omega\to0$ with ${\rm Im}\omega>0$ of 
$$	\int_0^\infty e^{i\omega t}dt
	\int\rho_a(dx)\,X(x)\cdot\nabla_x(A\circ f_a^t)      $$
\par}
\vfill\eject
	{\bf 1. Introduction.}
\medskip
	Given a time evolution $(x,t)\mapsto f^tx$, with $x\in\hbox{manifold }M$, $t\in{\bf R}$, it is often possible to find a set $S\subset M$ and an invariant probability measure $\rho$ on $M$ such that ${\rm lebesgue}(S)>0$ ({\it i.e.}, $S$ has positive Lebesgue measure), and 
$$	\lim_{T\to\infty}{1\over T}\int_0^TA(f^tx)\,dt=\rho(A)\qquad
	\hbox{if }x\in S\eqno{(1)}      $$
whenever $A:M\to{\bf R}$ is continuous.  Such measures $\rho$ are called SRB measures or SRB states.  (In the case of a discrete time dynamical system, the integral in $(1)$ is replaced by a sum).  
\medskip
	SRB measures were defined and studied by Ya. Sinai, D. Ruelle and R. Bowen for uniformly hyperbolic\footnote{*}{We call uniformly hyperbolic the Anosov systems [1] and the more general Axiom A systems introduced by Smale [32] (see also Bowen [7]).} systems [31], [24], [8].  Then the concept was extended to general smooth dynamical system by F. Ledrappier, J.-M. Strelcyn and L.-S. Young [18], [19].  Later it was found that, in a number of situations where specific geometric information is available, one can prove detailed properties of SRB measures (see in particular L.-S. Young [33], and the monograph by C. Bonatti, L. Diaz and M. Viana [3]).  
\medskip
	The SRB measures describe the statistical properties of physical systems, in particular in nonequilibrium statistical mechanics [28].  It is therefore desirable to study how these measures depend on parameters ({\it i.e.}, on the dynamical system $(f^t)$).  For the large systems of statistical mechanics, a {\it linear response} is often observed experimentally when parameters are varied.  This means that the expectation value $\rho(A)$ of an {\it observable} $A$ should depend differentiably on parameters.  It is not clear at present how to reconcile the concept of linear response with the fact that typical dynamical systems depend very discontinuously on parameters (and may exhibit a dense set of bifurcations).  The uniformly hyperbolic case is however amenable to discussion (in physical situations, this amounts to accepting the {\it chaotic hypothesis} of G. Gallavotti and E.G.D. Cohen [16]).  A formula for the derivative of SRB states with respect to parameters has been obtained in the case of Axiom A diffeomorphisms in [27].  Here we shall study Axiom A flows.  
\medskip
	A precise statement of our results is given as Theorem A and Theorem B below.  The general idea of the proofs is to use the symbolic dynamics for hyperbolic flows to study their SRB states, also applying methods of the thermodynamic formalism\footnote{**}{Ya. Sinai introduced Markov partitions, symbolic dynamics, and studied the ergodic theory for Anosov diffeomorphisms [29], [30], [31].  A partial generalization to flows was given by M. Ratner [23].  Then R. Bowen gave a general definition of Markov partitions for Axiom A diffeomorphisms [4] and flows [5].  The ergodic theory for Axiom A flows was studied by R. Bowen and D. Ruelle [8], introducing what are here called SRB states on attractors for Axiom A flows.  Some abstract results applicable to SRB states originate from equilibrium statistical mechanics and are subsumed in the so-called thermodynamic formalism [6], [25].}.  
\medskip
	It will be convenient to use the following notation for the derivative at $x$ of a function $A$ on the manifold $M$ in the direction of the vector field $X$: 
$$	X(x)\cdot\nabla_xA=(D_xA)X(x)      $$
If $f$ is a diffeomorphism of $M$ we have thus 
$$	X(x)\cdot\nabla_x(A\circ f)=(D_{fx}A)(T_xf)X(x)      $$
\medskip
	{\bf 2. Differentiability of SRB states for hyperbolic systems.}
\medskip
	Let $r\ge3$, and $(f_a^t)$ be a ${\cal C}^r$ hyperbolic dynamical system (diffeomorphism or flow) depending smoothly on a parameter $a$, with an SRB measure $\rho_a$.  There are a number of results on the smoothness of $a\mapsto\rho_a$ as a distribution, {\it i.e.}, of $a\mapsto\rho_a(A)$ when $A$ is smooth.  See [21], [17], [10], [11], [2].
\medskip
	For applications to statistical physics it is desirable to have an explicit expression for $d\rho_a(A)/da$.  In the case of an Axiom A diffeomorphism $f_a$, writing $X_a=({d\over da}f_a)\circ f_a^{-1}$, we obtain by a formal calculation
$$	{d\over da}\rho_a(A)
=\sum_{k=0}^\infty\int\rho_a(dx)\,X_a(x)\cdot\nabla_x(A\circ f_a^k)      $$
If $f_a$ is mixing, this result holds with an exponentially convergent sum over $k$, as shown in [27].  The proof is more difficult than one might anticipate.  (For other differentiability results see [14]).  
\medskip
	In the present paper we tackle the case of an Axiom A flow $(f_a^t)$ defined by a vector field ${\cal X}+aX$.  Here a formal calculation yields 
$$	{d\over da}\rho_a(A)
=\int_0^\infty dt\int\rho_a(dx)\,X(x)\cdot\nabla_x(A\circ f_a^t)      $$
What we shall show is that the Fourier transform 
	$$\int_0^\infty e^{i\omega t}dt
	\int\rho_a(dx)\,X(x)\cdot\nabla_x(A\circ f_a^t)      $$
(defined as a distribution) extends to a holomorphic function of $\omega$ near $\omega=0$ such that its value at 0 is ${d\over da}\rho_a(A)$.  
\medskip
	While the proofs presented here are relatively straightforward, they make detailed use of the references [5], [8], [25], [26], and lead to somewhat heavy formulas.  (The author has tried without success to find simpler and more direct arguments).  
\medskip
	{\bf 3. Theorem A.} 
\medskip
	{\sl Let ${\cal X}$ and $X$ be ${\cal C}^r$ vector fields ($r\ge3$) on the compact manifold $M$, and let $(f_a^t)$ be the flow defined by ${\cal X}+aX$.  We assume that for small $a$ the flow $(f_a^t)$ has a nontrivial\footnote{*}{\rm The attractor $\Lambda_a$ is nontrivial if it is not a fixed point or a periodic orbit.} Axiom A attractor $\Lambda_a$ (depending continuously on $a$) with SRB measure $\rho_a$.  
\medskip
	If $A\in{\cal C}^{r-1}(M)$, the function $a\mapsto\rho_a(A)$ is ${\cal C}^{r-2}$ and ${d\over da}\rho_a(A)|_{a=0}$ is the value at $\omega=0$ of the function defined for ${\rm Im}\omega>0$ by 
$$	\omega\mapsto\int_0^\infty e^{i\omega t}dt
	\int\rho_0(dx)\,X(x)\cdot\nabla_x(A\circ f_0^t)      $$
which extends meromorphically to $\{\omega:{\rm Im}\omega>-\delta\}$ for some $\delta>0$, without pole at $\omega=0$.}
\medskip
	Note that the theorem does not assume the flow $(f_a^t)$ to be mixing.  If $\int_0^\infty dt\,|\rho_0((A\circ f_0^t).C)|<\infty$, where $C={\rm div}_v^{cu}(X^c+X^u)$ is defined in Section 7 below, we have 
$$	{d\over da}\rho_a(A)|_{a=0}
=\int_0^\infty dt\int\rho_0(dx)\,X(x)\cdot\nabla_x(A\circ f_0^t)      $$
[There are a number of results on decay of correlations for hyperbolic flows, see in particular Chernov [9], Dolgopyat [12], [13], Liverani [20], Fields, Melbourne and T\"or\"ok [15].  Since $C$ is H\"older but not smooth in general, only [20] applies directly in the present situation].  
\medskip
	A proof of Theorem A will be obtained from Theorem B below.
\medskip
	{\bf 4. Corollary.} 
\medskip
	{\sl Suppose that the vector field $X_t$ is constant in $t$ and equal to $X$ when $t\le t_0$ for some time $t_0$, but that $X_t$ may depend (smoothly) on $t$ for $t\ge t_0$.  Write $f_a^{(t,t_0)}x_0=x(t)$ where ${d\over dx}x(t)={\cal X}(x(t))+aX_t(x(t))$ and $x(t_0)=x_0$.  One can then define a time dependent SRB state $\rho_a^t=f_a^{(t,t_0)}\rho_a$ so that it reduces to $\rho_a$ for $t\le t_0$.  With this definition, if $\int_0^\infty dt\,|\rho_0((A\circ f_0^t).C)|<\infty$, 
$$	{d\over da}\rho_a^t(A)|_{a=0}=\int_{-\infty}^td\tau
\int\rho_0(dx)\,X_{\tau}(x)\cdot\nabla_x(A\circ f_0^{t-\tau})      $$}
\indent
	The Corollary follows  directly from Theorem A when $t<t_0$.  To obtain the general case differentiate both sides with respect to $t$.

	Before we formulate Theorem B, we need some facts and definitions.
\medskip
	{\bf 5. Correlation functions.}  
\medskip
	If $B,B'$ are smooth functions on a neighborhood of $\Lambda_0$ in $M$, their correlation function is $t\mapsto\rho_{BB'}(t)=\rho_0((B\circ f_0^t).B')-\rho_0(B)\rho_0(B')$.  Multiplying by the characteristic function $\chi^+$ of $[0,+\infty)$ we obtain $\rho_{BB'}^+(t)=\rho_{BB'}(t)\chi^+(t)$, and taking the Fourier transform
$$	\hat\rho_{BB'}^+(\omega)=\int_0^\infty e^{i\omega t}dt
	[\rho_0((B\circ f_0^t).B')-\rho_0(B)\rho_0(B')]      $$
This is a distribution, boundary value of a holomorphic function in the upper half complex plane, which furthermore extends to a meromorphic function in $\{\omega:{\rm Im}\omega>-\delta'\}$ for some $\delta'>0$, with no pole at $\omega=0$, as discussed in [22], [26].  Actually, the discussion in [26] uses a {\it symbolic} representation of $\Lambda$: points have a description $(\xi,t)$ where $\xi$ belongs to a Cantor set $\Sigma$, and $t$ to an interval of ${\bf R}$.  Instead of smooth $B,B'$ one takes $B,B'\in{\cal C}^\sharp$, where ${\cal C}^\sharp$ is a Banach space of functions $t\mapsto B(\cdot,t)$, continuous: interval of ${\bf R}\to{\cal C}^\alpha(\Sigma)$.  

\noindent
[To make the connection with the formalism of [26], it is useful to know that if $t\mapsto B(\cdot,t),\zeta(\cdot,t)$ are continuous: interval $\to{\cal C}^\alpha(\Sigma)$, and $t\mapsto B(\cdot,t)$ is ${\cal C}^2$: interval $\to$ bounded functions on $\Sigma$, then $t\mapsto B(\cdot,\zeta(\cdot,t))$ is continuous: interval $\to{\cal C}^\alpha(\Sigma)$].  
\medskip
	For our purposes the function $B'=C$ to be introduced below will belong to ${\cal C}^\sharp$ rather than being smooth.  
\medskip
	{\bf 6. The volume elements $\tilde v$ and $v$.}
\medskip
	Let ${\cal V}^u$ denote a strong unstable manifold for the flow $(f_0^t)$.  We have thus ${\cal V}^u\subset\Lambda_0$, and ${\cal V}^u$ is $u$-dimensional.  There is a natural volume element $\tilde v$ on each such ${\cal V}^u$ so that, for all $t$, the natural volume element on $f_0^t{\cal V}^u$ is the image by $f_0^t$ of the measure $\tilde v$, up to a multiplicative constant.  This is seen in the same way as for the existence of a natural volume element on unstable manifolds contained in an attractor for an Axiom A diffeomorphism (see [27]).  Here again $\tilde v$ has ${\cal C}^{r-1}$ density, and is uniquely defined up to a multiplicative constant.  
\medskip
	If $\tilde{\cal V}^u$ is a $u$-dimensional manifold contained in a center-unstable manifold, and trans\-versal to the flow $(f_0^t)$, we can define a volume element $\tilde v$ on $\tilde{\cal V}^u$ as the image of $\tilde v$ on a strong unstable manifold ${\cal V}^u$ by a Poincar\'e map.  In this manner we obtain a natural volume element $\tilde v$, defined up to a multiplicative constant and corresponding to Poincar\'e maps acting on manifolds $\tilde{\cal V}^u$ transversal to $(f_0^t)$.  
\medskip
	Let now ${\cal W}^{cu}$ denote a center-unstable manifold for the flow $(f_0^t)$.  We have thus ${\cal W}^{cu}\subset\Lambda_0$, and ${\cal W}^{cu}$ is $(u+1)$-dimensional.  Take a chart $S\times I$ of $M$ such that ${\cal X}$ is the unit vector in the last coordinate direction, and $I$ is an interval of ${\bf R}$.  Assuming also $\tilde{\cal V}^u\subset S$ we may write locally ${\cal W}^{cu}=\tilde{\cal V}^u\times I$ and define 
$$	v=\tilde v\times{\rm lebesgue}      $$
A volume element $v$ is thus given on the center-unstable manifolds ${\cal W}^{cu}$, and is unique up to a multiplicative constant.  Note that $v$ has ${\cal C}^{r-1}$ density and that $f_0^t$ sends $v$ to $v$ up to a multiplicative constant.  [We shall see in Section 10 that $v$ is (up to a multiplicative constant) the conditional probability of the SRB measure $\rho_0$ on the (local) center-unstable manifold ${\cal W}^{cu}$].
\medskip
	{\bf 7.  The function $C={\rm div}_v^{cu}(X^c+X^u)$.}
\medskip
	For $x\in\Lambda_0$, let $T_xM=E_x^c+E_x^s+E_x^u$, where $E_x^c$ is $1$-dimensional containing ${\cal X}(x)$, and $E_x^s,E_x^u$ are the strong stable and unstable subspaces at $x$ for $(f_0^t)$.  We write 
$$	X(x)=X^c(x)+X^s(x)+X^u(x)      $$
with $X^c(x)\in E_x^c$, $X^s(x)\in E_x^s$, $X^u(x)\in E_x^u$.  If we take again a chart $S\times I$ of $M$ such that ${\cal X}$ is that unit vector in the last coordinate direction, we see that $E_x^c$ is independent of $x$, while $E_x^s,E_x^u$ depend H\"older continuously on $x$, and are independent of the last coordinate of $x$.  In particular $X^c(x),X^s(x),X^u(x)$ have ${\cal C}^r$ dependence on the last coordinate of $x$ (while depending H\"older continuously on $x$).  
\medskip
	The divergence of $X^c+X^u$ with respect to the volume element $v$ on the manifold ${\cal W}^{cu}$ is denoted by ${\rm div}_v^{cu}(X^c+X^u)$.  It is, {\it a priori}, a distribution, but we shall show that it is actually a H\"older continuous function on $\Lambda_0$ (note that this is a local question).  
\medskip
	Let $f_0^tx\in{\cal W}^{cu}$, with $x\in S\cap{\cal W}^{cu}=\tilde{\cal V}^u$.  We may write $X^c+X^u=X'^c+X'^u$ where $X'^c(f_0^tx)\in E_x^c$ and $X'^u(f_0^tx)\in T_xS\cap(E_x^c+E_x^u)$.  We have then ${\rm div}_v^{cu}(X^c+X^u)=\partial X'^c+{\rm div}_{\tilde v}X'^u$ where $\partial X'^c$ denotes the derivative of $X'^c$ with respect to the last coordinate ({\it i.e.}, $(\partial X)(f_0^tx)=\partial_tX(f_0^tx)$).  Since $\partial X$ is ${\cal C}^{r-1}$, $\partial X'^c$ is H\"older continuous.  Note that we may also write $X=X''^c+X''^s+X'^u$ where $X''^c(f_0^tx)\in E_x^c$ and $X''^s(f_0^tx)\in T_xS\cap(E_x^c+E_x^s)$.  The definition of ${\rm div}_{\tilde v}$ in ${\cal W}^{cu}\cap S$ is now very similar to that of ${\rm div}^u$ for the case of hyperbolic diffeomorphisms in [27], provided we replace the diffeomorphism $f$ by Poincare maps of $(f_0^t)$.  In fact, using a Markov partition for $(f_0^t)$ we see that we need only a finite number of Poincare maps $f_0^{T_{k\ell}}$ between sections $S_k$, $S_\ell$.  The stable and unstable directions for the system of Poincar\'e maps are $T_xS\cap(E_x^c+E_x^s)$, $T_xS\cap(E_x^c+E_x^u)$.  Using the results of [27] we obtain thus that ${\rm div}_{\tilde v}X'^u$ is H\"older, and therefore ${\rm div}_v^{cu}(X^c+X^u)$ is a H\"older function on $\Lambda_0$.  [Integration by parts will show (in Section 10) that $\rho_0(C)=0$ because boundary terms cancel out].  Instead of $X$ we may use $\partial X$ in the above argument, and find that 
$$	f_0^t x\mapsto\partial_tC(f_0^tx)
	={\rm div}_v^{cu}(\partial X^c+\partial X^u)(f_0^tx)      $$
is H\"older continuous on $\Lambda_0$.  From this results that $t\mapsto(x\mapsto C(f_0^tx))$ defines a ${\cal C}^1$ function to ${\cal C}^\alpha(S)$.
\medskip
	{\bf 8. Theorem B.} 
\medskip
	{\sl Under the conditions of Theorem A we have 
$$	{d\over da}\rho_a(A)|_{a=0}=
	\int\rho_0(dx)\,(D_xA)\int_0^\infty dt\,(T_{f_0^{-t}x}f_0^t)X^s(f_0^{-t}x)
	-\hat\rho_{AC}^+(0)      $$}
[If $\int_0^\infty dt\,|\rho_0((A\circ f_0^t).C)|<\infty$, we have $\hat\rho_{AC}^+(0)=\int_0^\infty\rho_0((A\circ f_0^t).C)$]
\medskip
	The proof of Theorem B will occupy most of the rest of this paper.  It is based on the study of SRB states with help of a Markov partition.  We start with the unperturbed dynamics ({\it i.e.}, $a=0$, the index $a$ will be omitted until Section 11).
\medskip
	Let thus, for $r\ge3$, $\Lambda$ be an Axiom A attractor for the flow $(f^t)$ defined on the manifold $M$ by the C$^r$ vector field ${\cal X}$: 
$$	{df^tx\over dt}={\cal X}(f^tx)\eqno{(2)}      $$
with $f^0x=x$.  There is a unique SRB measure $\rho$ with support $\Lambda$ for the flow $(f^t)$.  A perturbation $\delta{\cal X}$ of the vector field ${\cal X}$ causes a change $\delta\rho$ of the SRB state $\rho$ and we have formally 
$$	\delta\rho(A)=\int_0^\infty ds
\int\rho(dx)\delta{\cal X}(x)\cdot\nabla_x(A\circ f^s)\eqno{(3)}      $$
for smooth $A:M\to{\bf R}$.  The main purpose of the present paper is to provide a proof of a modified version of $(3)$, as described in Theorem A and Theorem B above.  
\medskip
	{\bf 9. Markov partition for the flow $(f^t)$.}  
\medskip
	We introduce a Markov partition with data as follows (see [5]).  A finite index set $J$ is given, and an $J\times J$ matrix $\tau$ with entries $0$ or $1$ such that all entries of some power of $\tau$ are $>0$.  We denote by $(\Sigma,\sigma)$ the mixing subshift of finite type defined by $J,\tau$, and let 
$$	\Sigma_k=\{(\xi_j)_{j\in{\bf Z}}:\xi_0=k\}\qquad,\qquad
	\Sigma_{k\ell}=\{(\xi_j)_{j\in{\bf Z}}:\xi_0=k,\xi_1=\ell\}      $$
The construction of the Markov partition uses small pieces $S_k$ of manifolds transversal to the flow $(f^t)$ for $k\in J$ (the $S_k$ are open codimension 1 smooth submanifolds of $M$).  When $\tau_{k\ell}=1$, an open subset $S_{k\ell}$ of $S_k$ and a C$^r$ real function $T_{k\ell}>0$ on $S_{k\ell}$ are given such that $f^{T_{k\ell}}S_{k\ell}\subset S_\ell$.  Finally, for some standard metric on $\Sigma$, there is a $\alpha$-H\"older continuous map $\pi:\Sigma\to\cup_k(S_k\cap\Lambda)$ such that 
\def\fl#1{\big\downarrow\vbox to 5mm{}\rlap{$\scriptstyle #1$}}
$$\matrix{\Sigma_{k\ell}&\buildrel\sigma\over{\longrightarrow}&\Sigma_\ell\cr
\fl{\pi}&&\fl{\pi}\cr
S_{k\ell}&\buildrel f^{T_{k\ell}}\over{\longrightarrow}&S_\ell\cr}$$
is commutative.  A positive $\alpha$-H\"older continuous function $\psi:\Sigma\to{\bf R}$ is defined by 
$$	\psi(\xi)=T_{k\ell}(\pi\xi)\qquad\hbox{when}\qquad
	\xi\in\Sigma_{k\ell}      $$
Also, if $A$ is H\"older continuous on $\Lambda$ we define a $\gamma$-H\"older continuous function $\tilde A$ on $\Sigma$ by 
$$	\tilde A(\xi)=\int_0^{\psi(\xi)}dt\,A(f^t\pi\xi)\eqno{(4)}      $$
(here $\gamma=\alpha$ if $A\in$C$^1(M)$, otherwise we have to choose some $\gamma\le\alpha$).
\medskip
	{\bf 10. Equilibrium states.}
\medskip
	We use here the formalism of [8], calling {\it equilibrium states} the invariant probability measures described elsewhere as {\it Gibbs states}.  The {\it pressure} of a H\"older continuous function $\phi:\Lambda\to{\bf R}$ with respect to the flow $(f^t)$ is 
$$	c=\sup_\nu{h_\sigma(\nu)+\nu(\tilde\phi)\over\nu(\psi)}      $$
where the sup is over $\sigma$-invariant probability measures $\nu$ on $\Sigma$, $h_\sigma$ denotes the {\it entropy} with respect to the shift $\sigma$, and $\tilde\phi$ is defined according to (4).  Let $\nu_0$ be the unique equilibrium state for $\tilde\phi-c\psi$ on $\Sigma$.  Then the unique equilibrium state $\mu_\phi$ of $\phi$ for the flow $(f^t)$ on $\Lambda$ is given by 
$$	\mu_\phi(A)={\nu_0(\tilde A)\over\nu_0(\psi)}\eqno{(5)}      $$
We shall be interested in the case when $\phi=\phi^{(u)}$ is minus the time derivative of the unstable Jacobian: 
$$	\phi=\phi^{(u)}
=-{d\over dt}\lambda_t^+|_{t=0}=-{d\over dt}\log\lambda_t^+|_{t=0}      $$
with 
$$	\lambda_t^+(x)
	=||(T_xf^t)^{\wedge(u+1)}|\hbox{volume element of ${\cal W}^{cu}$}||
	=||(T_xf^t)^{\wedge u}|\hbox{volume element of ${\cal V}^u$}||      $$
Notice that we have 
$$	\phi^{(u)}(f^tx)=-{d\over dt}\log\lambda_t^+(x)      $$
\indent
	For $\phi^{(u)}$ one can show that the pressure vanishes ($c=0$) and $\mu_{\phi^{(u)}}$ is the SRB measure $\rho$ on $\Lambda$ for $(f^t)$.  Details and proofs of the above construction of the SRB measure $\rho$ are given in [8].  Note that the function $\tilde\phi$ corresponding to $\phi=\phi^{(u)}$ is -- up to a minus sign and composition with $\pi$ -- the unstable Jacobian $(\lambda_{T_{k\ell}}^+)$ of $(f^{T_{k\ell}})$ acting on $(S_\ell)$.  This reduces the study of $\nu_0$ to the situation discussed in [27] for an Axiom A diffeomorphism $f$, with the replacement of $f$ by $(f^{T_{k\ell}})$.  In particular (5) shows that the conditional measures of $\rho$ on ${\cal W}^{cu}$ are of the form $v=\tilde v\times$ lebesgue.  We obtain thus $\rho(C)=\rho({\rm div}_v^{cu}(X^c+X^u))=0$ by integration by parts and cancellation of boundary terms, as announced earlier.
\medskip
	{\bf 11. Flows depending on a parameter $a$.}
\medskip
	If we replace ${\cal X}$ in $(2)$ by ${\cal X}+aX$ for $a\in(-\epsilon,\epsilon)$ we may leave $\Sigma$, $\sigma$, $S_k$, $S_{k\ell}$ unchanged but replace $(f^t)$, $\Lambda$, $T_{k\ell}$, $\pi$, $\psi$, $\phi$, $\tilde A$ by $(f_a^t)$, $\Lambda_a$, $T_{ak\ell}$, $\pi_a$, $\psi_a$, $\phi_a$, $\tilde A_a$.  Call $\pi_*$ the map $\pi$ introduced in Section 9.  A hyperbolic fixed point argument shows that for suitable $\alpha>0$ there is a $\alpha$-H\"older $\pi_a:\Sigma\to\cup S_\ell$ such that 
$$	f_a^{T_{ak\ell}}\circ\pi_a\circ\sigma^{-1}=\pi_a
	\qquad\hbox{on}\qquad\sigma\Sigma_{k\ell}      $$
and $a\mapsto\pi_a$ is ${\cal C}^{r-1}:(-\epsilon,\epsilon)\to{\cal C}^\alpha(\Sigma\to\cup S_\ell)$, reducing to $\pi_*$ for $a=0$.
\medskip
	Here are details.  Define $\Psi_a=(\Psi_{ak\ell})$ where 
$$	\Psi_{akl}\pi=f_a^{T_{ak\ell}}\circ\pi\circ\sigma^{-1}
	\qquad\hbox{on}\qquad\sigma\Sigma_{k\ell}      $$
for $(a,\pi)$ close to $(0,\pi_*)$.  Then $\Psi_a$ maps a neighborhood of $\pi_*$ in the H\"older space ${\cal C}^\alpha(\Sigma\to\cup_{k\ell}S_{k\ell})$ to ${\cal C}^\alpha(\Sigma\to\cup_{k\ell}S_{k\ell})$.  We assume that we have charts identifying the $S_{k\ell}$ with open subsets of ${\bf R}^{{\rm dim}M-1}$, so that ${\cal C}^\alpha(\Sigma\to\cup_{k\ell}S_{k\ell})\subset$ ${\cal C}^\alpha(\Sigma\to{\bf R}^{{\rm dim}M-1})$.  Note that $(a,\pi)\mapsto\Psi_a\pi$ is ${\cal C}^{r-1}$ hence ${\cal C}^1$ from a neighborhood of $(0,\pi_*)$ in ${\bf R}\times{\cal C}^\alpha(\Sigma\to{\bf R}^{{\rm dim}M-1})$ to ${\cal C}^\alpha(\Sigma\to{\bf R}^{{\rm dim}M-1})$.  Taking $a=0$ we see that $\pi_*$ is a fixed point of $\Psi_0$ (see the commutative diagram in Section 9 above).  The derivative $D_{\pi_*}\Psi_0$ is a bounded linear operator on ${\cal C}^\alpha(\Sigma\to{\bf R}^{{\rm dim}M-1})$.  Let $V_{\pi_*\xi}^s,V_{\pi_*\xi}^u\subset{\bf R}^{{\rm dim}M-1}$ denote the stable and unstable subspaces at $\pi_*\xi$.  (When $\xi\in\Sigma_\ell$ these are the intersections with $T_{\pi_*\xi}S_\ell$ of the  center-stable and center-unstable spaces at $\pi_*\xi$ for $(f_0^t)$, or the stable and unstable spaces for the $f_0^{T_{0k\ell}}$ acting on $\cup_\ell S_\ell$).  We have chosen $\alpha>0$ such that $\pi_*$ is $\alpha$-H\"older, and we may assume that also $\xi\mapsto V_{\pi_*\xi}^{s,u}$ is $\alpha$-H\"older.  The spaces $V_*^{s,u}$, defined to consist of the $\alpha$-H\"older maps $\xi\to V_{\pi_*\xi}^{s,u}$ are closed linear subspaces of ${\cal C}^\alpha(\Sigma\to{\bf R}^{{\rm dim}M-1})$, and ${\cal C}^\alpha(\Sigma\to{\bf R}^{{\rm dim}M-1})=V_*^s\oplus V_*^u$.  
\medskip
	We show now that $D_{\pi_*}\Psi_0$ is a hyperbolic operator with respect to the direct sum decomposition $V_*^s\oplus V_*^u$, provided $\alpha$ has been chosen small enough, {\it i.e.}, if $\alpha$ is replaced by a suitable $\beta$ (with $0<\beta<\alpha$) which we shall now determine.  It suffices to prove that $D_{\pi_*}\Psi_0$ induces a contraction on $V_*^s$, where $D_{\pi_*}\Psi_0$ is the map 
$$	u\mapsto(Tf_0^{T_{0k\ell}})(u\circ\sigma^{-1})      $$
Using an ``adapted metric'' on $M$ we may assume for the uniform norm 
$$	||Tf_0^{T_{0k\ell}}|\hbox{stable direction}||_0\le\lambda<1      $$
In the definition of the ${\cal C}^\beta$ norm 
$$	||\Phi||=\max\big(\sup_\xi|\Phi(\xi)|,
\sup_{\xi\ne\eta}{|\Phi(\xi)-\Phi(\eta)|\over d(\xi,\eta)^\beta}\big)      $$
we take the second $\sup$ only over pairs $(\xi,\eta)$ such that $d(\xi,\eta)^\beta)<\epsilon$, where the constant $\epsilon$ will be fixed later (small but $>0$).
\medskip
	Write $T_\xi=T_{\pi_*\xi}f_0^{T_{0k\ell}}$, $\delta=d(\xi,\eta)$.  Given $u\in V_*^s$ (with ${\cal C}^\beta$ norm $||u||$) we may for each pair $(\xi,\eta)$ with small $\delta$ choose $v\in V_{\pi_*\xi}^s$ with $|v-u(\eta)|\le||u||O(\delta^\alpha)$.  We have
$$	T_\xi u(\xi)-T_\eta u(\eta)
	=T_\xi(u(\xi)-v)+T_\xi v-T_\eta v+T_\eta(v-u(\eta))      $$
$$	|T_\xi(u(\xi)-v)|
\le\lambda|u(\xi)-v|\le\lambda|u(\xi)-u(\eta)|+||u||O(\delta^\alpha)      $$
$$	|T_\xi v-T_\eta v|\le||u||O(\delta^\alpha)      $$
$$	|T_\eta(v-u(\eta))|\le||u||O(\delta^\alpha)      $$
hence
$$	|T_\xi u(\xi)-T_\eta u(\eta)|
	\le||u||(\lambda\delta^\beta+O(\delta^\alpha))      $$
Since $d(\sigma\xi,\sigma\eta)\ge C\delta$ we have
$$	{|T_\xi u(\xi)-T_\eta u(\eta)|\over d(\sigma\xi,\sigma\eta)^\beta}
	\le||u||{\lambda\delta^\beta+O(\delta^\alpha)\over C^\beta\delta^\beta}
	=||u||({\lambda\over C^\beta}+O(\delta^{\alpha-\beta}))      $$
For small $\beta$ we have $\lambda/C^\beta<1$, and we may take $\epsilon$ such that 
$$	\lambda/C^\beta+O(\delta^{\alpha-\beta})<1
	\qquad{\rm if}\qquad0<\delta<\epsilon      $$
This concludes the proof that $D_{\pi_*}\Psi_0$ is hyperbolic for suitable $\beta$, {\it i.e.}, when $\alpha$ is chosen small enough.  We may thus apply the implicit function theorem to obtain the existence of $\pi_a$ with the properties indicated above.  
\medskip
	{\bf 12. Smooth dependence of SRB state with respect to $a$.}
\medskip
	Let $\phi_a=\phi_a^{(u)}$ be minus the time derivative of the unstable Jacobian for $(f_a^t)$ and $\nu_a$ the unique equilibrium state for $\tilde\phi_a$ on $\Sigma$, where 
$$	\tilde\phi_a(\xi)=\int_0^{\psi_a(\xi)}dt\,\phi_a(f_a^t\pi_a\xi)      $$
Then, according to Section 10, the SRB measure $\rho_a$ for $(f_a^t)$ on $\Lambda_a$ is given by 
$$	\rho_a(A)={\nu_a(\tilde A_a)\over\nu_a(\psi_a)}      $$
Assuming $A\in{\cal C}^r(M)$ we find that $a\mapsto\psi_a,\tilde A_a$ are ${\cal C}^{r-1}:(-\epsilon,\epsilon)\to{\cal C}^\alpha(\Sigma)$ because we know that $a\mapsto\pi_a$ is ${\cal C}^{r-1}$, and
$$	\psi_a(\xi)
	=T_{ak\ell}(\pi_a\xi)\qquad\hbox{for}\qquad\xi\in\Sigma_{k\ell}      $$
$$	\tilde A_a(\xi)=\int_0^{\psi_a(\xi)}dt\,A(f_a^t\pi_a\xi)      $$
\indent
	The set $\hat\Lambda_a=E_{\Lambda_a}^u$ of unstable subspaces is an Axiom A attractor for the ${\cal C}^{r-1}$ action of $(Tf_a^t)$ on the Grassmannian $\widehat M\to M$.  Therefore if $\hat\pi_a:\Sigma\to\hat\Lambda_a$ makes the diagram 
\def\fl#1{\downarrow\vbox to 5mm{}\rlap{$\scriptstyle #1$}}
$$\matrix{&&\hat\Lambda_a\cr
&\buildrel\hat\pi_a\over\nearrow&\fl{}\cr
\Sigma&\buildrel\pi_a\over{\rightarrow}&\Lambda_a\cr}$$
commutative, we see that $a\mapsto\hat\pi_a$ is ${\cal C}^{r-2}:(-\epsilon,\epsilon)\to{\cal C}^\alpha$ (where we may again have to replace the current value of $\alpha$ by a lower one).  Note that 
$$	\tilde\phi_a(\xi)=-\log\lambda_{\psi_a(\xi)}^+(\pi_a\xi)      $$
where $\lambda_t^+(\pi_a\xi)$ is the unstable Jacobian $||(T_{\pi_a\xi}f_a^t)^{\wedge u}|\hbox{volume element of }\hat\pi_a\xi||$.  Note that $\lambda_{\psi_a(\xi)}^+(\pi_a\xi)$ is a ${\cal C}^{r-1}$ function of $a,\psi_a(\xi),\hat\pi_a\xi$, hence $a\mapsto\tilde\phi_a(\cdot)$ is ${\cal C}^{r-2}:(-\epsilon,\epsilon)\to{\cal C}^\alpha(\Sigma\to{\bf R})$.  Therefore $a\mapsto\nu_a$ is ${\cal C}^{r-2}:(-\epsilon,\epsilon)\to({\cal C}^\alpha(\Sigma\to{\bf R}))^*$.  [We use here the thermodynamic formalism to obtain the ${\cal C}^\omega$ dependence of $\nu_a$ (considered as an element of the Banach space dual of ${\cal C}^\alpha$) on $\tilde\phi_\alpha$ (considered as an element of ${\cal C}^\alpha$), see [25], Theorem 5.26].  Thus if $A\in{\cal C}^{r-1}(M)$, the function $a\mapsto\rho_a(A)=\nu_a(\tilde A_a)/\nu_a(\psi_a)$ is ${\cal C}^{r-2}$.
\medskip
	{\bf 13. Differentiating $a\mapsto\rho_a(A)$ at $a=0$.}
\medskip
	Writing $B=A-\rho_0(A)$ we have 
$$	\rho_a(A)=\rho_0(A)+\rho_a(B)
	=\rho_0(A)+{\nu_a(\tilde B_a)\over\nu_a(\psi_a)}      $$
where $\tilde B_a=\tilde A_a-\rho_0(A)\psi_a$.  Therefore 
$$	{d\over da}\rho_a(A)|_{a=0}
	={1\over\nu_0(\psi_0)}{d\over da}(\nu_a(\tilde B_a))|_{a=0}      $$
because $\nu_0(\tilde B_0)=\nu_0(\psi_0).\rho_0(B)=0$.  In view of the above formula we shall now study $\nu_a(\tilde B_a)$ to first order in $a$.  
\medskip
	{\bf 14. Reparametrization: modifying the map $\pi_a$ to first order in $a$.}
\medskip 
	A Markov partition parametrizes points of $\Lambda$ in the form $f^t\pi\xi$ where $\xi\in\Sigma$ and $0\le t<\psi(\xi)$.  We have taken $\pi\xi$ in a piece of smooth manifold $S_k$ transversal to the flow.  But we may just as well use a parametrization $f^t\pi^\sharp\xi$ of $\Lambda$, where $\pi^\sharp\xi=f^{\tau(\xi)}\pi\xi$ with continuous $\tau:\Sigma\to{\bf R}$.  
\medskip
	We consider a first such reparametrization which consists in replacing $S_k$ by a union of strong unstable manifolds (as is needed for the application of [26]).  This reparametrization corresponds to a H\"older continuous choice of $\xi\mapsto\tau(\xi)$, and replaces the $S_k$ by non-smooth ``manifolds'' in general.
\medskip
	We return now to smooth $S_k$ and write 
$$	\pi_a\xi=\pi_0\xi+a(U^c(\xi)+U^s(\xi)+U^u(\xi))      $$
to first order in $a$, with $U^c(\xi)\in E_{\pi_0\xi}^c$, $U^s(\xi)\in E_{\pi_0\xi}^s$, $U^u(\xi)\in E_{\pi_0\xi}^u$.  We may thus consider a second reparametrization:
$$	\pi_a^\sharp\xi=\pi_0\xi+a(U^s(\xi)+U^u(\xi))      $$
$$	=\pi_a\xi-aU^c(\xi)=f_a^{-a\theta(\xi)}\pi_a\xi      $$
where $\theta$ is defined by $U^c(\xi)=\theta(\xi){\cal X}(\pi_0\xi)$.  Note that the replacement of $\pi_a$ by $\pi_a^\sharp$ replaces also $\psi_a(\xi)$ by $\psi_a(\xi)+a\theta(\xi)-a\theta(\sigma\xi)$, $\tilde A_a(\xi)$ by $\tilde A_a(\xi)+a\theta(\xi)A(\pi_a\xi)-a\theta(\sigma\xi)A(\pi_a\sigma\xi)$, and $\tilde\phi_a(\xi)$ by $\tilde\phi_a(\xi)+a\theta(\xi)\phi_a(\pi_a\xi)-a\theta(\sigma\xi)\phi_a(\pi_a\sigma\xi)$.  Thus, the replacement of $\pi_a$ by $\pi_a^\sharp$ changes $\psi_a$, $\tilde A_a$, $\tilde\phi_a$ by a coboundary.  In particular $\nu_a$ and $\nu_a(\tilde B_a)$ are unchanged.  
\medskip
	Let us now perform the first and then the second reparametrization, {\it i.e.}, first replacing $S_k$ by a union of strong stable manifolds, and second taking 
$$	\pi_a^\sharp\xi=\pi_0\xi+a(U^s(\xi)+U^u(\xi))      $$
Here we have 
$$	\pi_a^\sharp\xi=\pi_a\xi-U^c(a,\xi)=f_a^{-\theta(a,\xi)}\pi_a\xi     $$
but because of the lack of smoothness of $S_k$, we cannot write $U^c(a,\xi)=aU^c(\xi)$, $\theta(a,\xi)=a\theta(\xi)$ in general.  Nevertheless, the replacement of $\pi_a$ by $\pi_a^\sharp$ changes $\psi_a,\tilde A_a,\phi_a$ by a coboundary, so that $\nu_a$ and $\nu_a(\tilde B_a)$ are unchanged.  In view of this we shall from now on replace $\pi_a$ by $\pi_a^\sharp$ and change $\psi_a$, $\tilde A_a$, $\tilde\phi_a$ accordingly, but without altering the notation.  
\medskip
	{\bf 15. Calculation of $\tilde B_a-\tilde B_0$.}
\medskip
	We have 
$$	\tilde B_a(\xi)-\tilde B_0(\xi)
	=\int_0^{\psi_a(\xi)}d\tau\,B(f_a^\tau(\pi_0\xi+aU^s(\xi)+aU^u(\xi))
	-\int_0^{\psi_0(\xi)}dt\,B(f_0^t\pi_0\xi)      $$
Write $X^c(x)=\eta(x){\cal X}(x)$, where $\eta$ is H\"older continuous on $\Lambda_0$ (and $\eta(f_0^t\pi_0\xi)$ is a smooth function of $t$).  We can then define a map $[0,\psi_a(\xi)]\to[0,\psi_0(\xi)]$ by $\tau\to t$ such that 
$$	{dt\over d\tau}=1+a\eta(f_0^\tau\pi_0\xi)      $$
Writing also $f_a^\tau=f_{a*}^t$ we obtain (to first order in $a$)
$$	\tilde B_a(\xi)-\tilde B_0(\xi)
	=\int_0^{\psi_0(\xi)}dt[(1-a\eta(f_0^t\pi_0\xi))B(f_{a*}^t(\pi_0\xi+aU^s(\xi)+aU^u(\xi))-B(f_0^t\pi_0\xi)]      $$
$$	=a(Z'-Z'')      $$
with 
$$	aZ'=\int_0^{\psi_0(\xi)}dt
	[B(f_{a*}^t(\pi_0\xi+aU^s(\xi)+aU^u(\xi))-B(f_0^t\pi_0\xi)]      $$
$$	Z''=\int_0^{\psi_0(\xi)}dt\,\eta(f_0^t\pi_0\xi)B(f_0^t\pi_0\xi)      $$
The contributions of $Z'$ and $Z''$ are evaluated in the Appendix.  
\medskip
	From now on we shall write $\pi$, $f$, $\psi$, $\nu$ instead of $\pi_0$, $f_0$, $\psi_0$, $\nu_0$.  For $n\ge0$, $\xi\in\Sigma$, we define
$$	\Psi(-n,\xi)=-\psi(\sigma^{-n}\xi)-\ldots-\psi(\sigma^{-1}\xi)      $$
$$	\Psi(n,\xi)=\psi(\xi)+\ldots+\psi(\sigma^{n-1}\xi)      $$
so that $\Psi(-n,\sigma^n\xi)=-\Psi(n,\xi)$, $\Psi(0,\xi)=0$, $\Psi(1,\xi)=\psi(\xi)$, and $f^{\Psi(k,\xi)}\pi\xi=\pi\sigma^k\xi$.  With this notation, the evaluation of $Z'$, $Z''$ in the Appendix yields the following result.
\medskip
	{\bf 16. Lemma.} {\sl We have }
$$	\nu({d\over da}\tilde B_a)|_{a=0}=\nu(Z'-Z'')      $$
$$	=\sum_{k=-\infty}^{-1}\int\nu(d\xi)\int_0^{\psi(\xi)}dt\,
	(D_{f^t\pi\xi}B)\int_{\Psi(k,\xi)}^{\Psi(k+1,\xi)}d\theta\, 
	(T_{f^\theta\pi\xi}f^{t-\theta})X^s(f^\theta\pi\xi)     $$
$$	+\int\nu(d\xi)\int_0^{\psi(\xi)}dt\,(D_{f^t\pi\xi}B)\int_0^t d\theta\, 
	(T_{f^\theta\pi\xi}f^{t-\theta})X^s(f^\theta\pi\xi)     $$
$$	-\sum_{k=-\infty}^{-1}\int\nu(d\xi)\int_0^{\psi(\xi)}dt\,B(f^t\pi\xi)
	\int_{\Psi(k,\xi)}^{\Psi(k+1,\xi)}d\theta\,
	({\rm div}_v^{cu}X^c)(f^\theta\pi\xi)      $$
$$	-\int\nu(d\xi)\int_0^{\psi(\xi)}dt\,B(f^t\pi\xi)\int_0^td\theta\,
	({\rm div}_v^{cu}X^c)(f^\theta\pi\xi)      $$
$$	-\int\nu(d\xi)\int_0^{\psi(\xi)}dt\,
	(D_{f^t\pi\xi}B)\int_t^{\psi(\xi)}d\theta\, 
	(T_{f^\theta\pi\xi}f^{t-\theta})X^u(f^\theta\pi\xi)      $$
$$	-\sum_{k=1}^\infty\int\nu(d\xi)\int_0^{\psi(\xi)}dt\,
	(D_{f^t\pi\xi}B)\int_{\Psi(k,\xi)}^{\Psi(k+1,\xi)}d\theta\, 
	(T_{f^\theta\pi\xi}f^{t-\theta})X^u(f^\theta\pi\xi)      $$
[The meaning of ${\rm div}_v^{cu}$ has been discussed in Section 7, the sums over $k$ converge exponentially, by hyperbolicity (directly) for the $X^s$ and $X^u$ parts, and by exponential decay of correlations for the $X^c$ part].  
\medskip
	{\bf 17. Evaluation of $\tilde\phi_a-\tilde\phi_0$.}
\medskip
	We have seen in Section 10 that the function $\tilde\phi$ corresponding to $\phi=\phi^{(u)}$ is -- up to a minus sign and composition with $\pi$ -- the unstable Jacobian $(\lambda_{T_{k\ell}}^+)$ of $(f^{T_{k\ell}})$ acting on $(S_\ell)$.  This reduces the study of $\nu$ to the situation discussed in [27] for an Axiom A diffeomorphism $f$, with the replacement of $f$ by $(f^{T_{k\ell}})$.  This remark remains true in the $a$-dependent situation, and reduces the evaluation of $\tilde\phi_a-\tilde\phi_0$ to the situation discussed in [27] for Axiom A diffeomorphisms.  We shall thus simply quote Proposition 1 of [27]II,which takes here the form 
$$	-{\tilde\phi_a-\tilde\phi_0\over\tilde\phi_0}
	\sim a({\rm div}_{\tilde v}^u\tilde X^u)\circ\pi      $$
In this formula the left-hand side is evaluated to first order in $a$, and we have used the following notation.  The equivalence $\sim$ means that the integrals of both sides with respect to every $\sigma$-invariant measure on $\Sigma$ coincide.  We have written 
$$	\int_0^{T_{k\ell}(x)}dt(T_{f^tx}f^{T_{k\ell}(x)-t})X^u(f^tx)
	=\tilde X^u(f^{T_{k\ell}(x)}x)      $$
Finally, the divergence ${\rm div}_{\tilde v}^u$ is computed, on the intersection ${\cal V}^u$ with $S_k$ of a center unstable manifold ${\cal W}^{cu}$, with respect to a natural volume element $\tilde v$ defined earlier.  (Note that, by our choice of $S_k$, ${\cal V}^u$ is a strong unstable manifold).  As in [27], and as in Section 7, ${\rm div}_{\tilde v}^u\tilde X^u$ is a H\"older continuous function on $S_k\cap\Lambda$. 
\medskip
	The relation between $\tilde X^u,X^u$ and $\tilde v,v$ also gives (see Section 7)
$$	({\rm div}_{\tilde v}^u\tilde X^u)(f^{T_{k\ell}(x)}(x))
	=\int_0^{T_{k\ell}(x)}dt\,({\rm div}_v^{cu}X^u)(f^tx)      $$
Therefore we may write 
$$	{d\over da}\log\tilde\phi_a(\xi)|_{a=0}\sim
-\int_0^{\psi(\xi)}dt\,({\rm div}_v^{cu}X^u)(f^t\pi\xi)=\gamma(\xi)      $$
The right-hand side is a H\"older continuous function of $\xi$ and, since $\nu_a$ is the equilibrium state for $\tilde\phi_a$, the thermodynamic formalism (see [25] Chapter 5, Exercise 5(b)) yields 
$$	{d\over da}\nu_a(\tilde B)|_{a=0}
	=\sum_{k=-\infty}^\infty[\nu(\tilde B.(\gamma\circ\sigma^k))
	-\nu(\tilde B)\nu(\gamma)]      $$
where the sum converges exponentially and, since $\nu(\tilde B)=0$, we find
$$	{d\over da}\nu_a(\tilde B)|_{a=0}
	=-\sum_{k=-\infty}^\infty\int\nu(d\xi)\tilde B(\xi)
\int_0^{\psi(\sigma^k\xi)}dt\,({\rm div}_v^{cu}X^u)(f^t\pi\sigma^k\xi)      $$
This yields the following result
\medskip
	{\bf 18. Lemma.} {\sl We have 
$$	{d\over da}\nu_a(\tilde B)|_{a=0}
	=-\sum_{k=-\infty}^\infty\int\nu(d\xi)\tilde B(\xi)\int_{\Psi(k,\xi)}
	^{\Psi(k+1,\xi)}d\theta\,({\rm div}_v^{cu}X^u)(f^\theta\pi\xi)      $$
$$  =-\sum_{k=-\infty}^\infty\int\nu(d\xi)\int_0^{\psi(\xi)}dt\,B(f^t\pi\xi)
	\int_{\Psi(k,\xi)}^{\Psi(k+1,\xi)}d\theta\,
	({\rm div}_v^{cu}X^u)(f^\theta\pi\xi)      $$
where the sum over $k$ converges exponentially.
\medskip
	The right-hand side above may be written as the sum of a part $Z_-$ where $\theta\le t$ and a part $Z_+$ where $\theta>t$.  In fact we claim that 
$$	{d\over da}\nu_a(\tilde B)|_{a=0}=Z_-+Z_+      $$
$$	Z_-=-\sum_{k=-\infty}^{-1}\int\nu(d\xi)\int_0^{\psi(\xi)}dt\,
	B(f^t\pi\xi)\int_{\Psi(k,\xi)}^{\Psi(k+1,\xi)}d\theta\,
	({\rm div}_v^{cu}X^u)(f^\theta\pi\xi)      $$
$$	-\int\nu(d\xi)\int_0^{\psi(\xi)}dt\,B(f^t\pi\xi)\int_0^td\theta\,
	({\rm div}_v^{cu}X^u)(f^\theta\pi\xi)      $$
$$	Z_+=\sum_{k=1}^\infty\int\nu(d\xi)\int_0^{\psi(\xi)}dt\,
	(D_{f^t\pi\xi}B)\int_{\Psi(k,\xi)}^{\Psi(k+1,\xi)}d\theta\,
	(T_{f^\theta\pi\xi}f^{t-\theta})X^u(f^\theta\pi\xi)      $$
$$	+\int\nu(d\xi)\int_0^{\psi(\xi)}dt\,(D_{f^t\pi\xi}B)
	\int_t^{\psi(\xi)}d\theta\,(T_{f^\theta\pi\xi}f^{t-\theta})
	X^u(f^\theta\pi\xi)   $$}
\indent
	For the calculation of the term $Z_+$, notice that if we write $\pi\xi=x$, the integral over $\nu(d\xi)\,dt$ reduces on the manifolds ${\cal W}^{cu}$ to integration over $\tilde v(dx)\,dt=v(dx\,dt)=dx_1\ldots dx_u\,dt$ for a suitable choice of coordinates.  Then, writing $X^u=Y$,
$$	B(f^tx)({\rm div}_v^{cu}X^u)(f^\theta x)
	=B(x,t)\sum_{k=1}^u\partial_kY^k(x,\theta)      $$
An integration by parts transforms this to 
$$	-\sum_{k=1}^u\partial_kB(x,t)Y^k(x,\theta)
	=-(D_{f^tx}B)(T_{f^\theta x}f^{t-\theta})X^u(f^\theta x)      $$
plus boundary terms involving $B(x,t)(T_{f^\theta x}f^{t-\theta})X^u(f^\theta x)$ with exponentially convergent integral over $\theta$.  The boundaries of pieces of ${\cal W}^{cu}$ are compact with zero measure, and it is readily seen that the boundary terms cancel.  
\medskip
	Putting Lemma 16 and Lemma 18 together yields: 
\medskip
	{\bf 19. Proposition.} 
\medskip
	{\sl We have 
$$	{d\over da}\nu_a(\tilde B_a)|_{a=0}
	=\sum_{-\infty}^{-1}\int\nu(d\xi)\int_0^{\psi(\xi)}dt\,
	(D_{f^t\pi\xi}B)\int_{\Psi(k,\xi)}^{\Psi(k+1,\xi)}d\theta\, 
	(T_{f^\theta\pi\xi}f^{t-\theta})X^s(f^\theta\pi\xi)     $$
$$	+\int\nu(d\xi)\int_0^{\psi(\xi)}dt\,(D_{f^t\pi\xi}B)\int_0^t d\theta\, 
	(T_{f^\theta\pi\xi}f^{t-\theta})X^s(f^\theta\pi\xi)     $$
$$	-\sum_{-\infty}^{-1}\int\nu(d\xi)\int_0^{\psi(\xi)}dt\,
B(f^t\pi\xi)\int_{\Psi(k,\xi)}^{\Psi(k+1,\xi)}d\theta\,C(f^\theta\pi\xi)   $$
$$	-\int\nu(d\xi)\int_0^{\psi(\xi)}dt\,
	B(f^t\pi\xi)\int_0^td\theta\,C(f^\theta\pi\xi)      $$
where we have written $C={\rm div}_v^{cu}(X^c+X^u)$}.
\medskip
	{\bf 20. Proof of Theorem A and Theorem B.}
\medskip
	We may write 
$$	\sum_{-\infty}^{-1}\int\nu(d\xi)\int_0^{\psi(\xi)}dt\,
	(D_{f^t\pi\xi}B)\int_{\Psi(k,\xi)}^{\Psi(k+1,\xi)}d\theta\, 
	(T_{f^\theta\pi\xi}f^{t-\theta})X^s(f^\theta\pi\xi)     $$
$$	+\int\nu(d\xi)\int_0^{\psi(\xi)}dt\,(D_{f^t\pi\xi}B)\int_0^t d\theta\, 
	(T_{f^\theta\pi\xi}f^{t-\theta})X^s(f^\theta\pi\xi)     $$
$$	=\int\nu(d\xi)\int_0^{\psi(\xi)}dt\,(D_{f^t\pi\xi}B)\int_{-\infty}^t 
d\theta\, (T_{f^\theta\pi\xi}f^{t-\theta})X^s(f^\theta\pi\xi)     $$
$$	=\nu(\psi)\int\rho(dx)\,(D_xB)
\int_0^\infty d\tau(T_{f^{-\tau}x}f^\tau)X^s(f^{-\tau}x)\eqno{(6)}      $$
This gives the first term occuring in theorem B.  In view of the exponential convergence of the integral over $\tau$, this term is also the value at $0$ of the expression 
$$	\omega\mapsto\int\rho(dx)\,(D_xB)\int_0^\infty e^{i\omega\tau}d\tau\,
	(T_{f^{-\tau}x}f^\tau)X^s(f^{-\tau}x)      $$
$$=\int_0^\infty e^{i\omega t}dt\int\rho(dx)\,X^s(x)\cdot\nabla_x(B\circ f^t)$$
which is holomorphic in $\omega$ for ${\rm Im}\omega>-\delta$, for some $\delta>0$, as required for Theorem A.
\medskip
	As to the series 
$$	-\sum_{-\infty}^{-1}\int\nu(d\xi)\int_0^{\psi(\xi)}dt\,
B(f^t\pi\xi)\int_{\Psi(k,\xi)}^{\Psi(k+1,\xi)}d\theta\,C(f^\theta\pi\xi) $$
$$	-\int\nu(d\xi)\int_0^{\psi(\xi)}dt\,B(f^t\pi\xi)
	\int_0^t d\theta\,C(f^\theta\pi\xi)\eqno{(7)}      $$
its sum is formally
$$	-\nu(\psi)\int_0^\infty dt\int\rho(dx)\,B(x)C(f^{-t}x)      $$
To obtain a rigorous estimate of (7) we consider the Fourier transform, as temperate distribution, of $\rho_{BC}^+(\cdot)=\rho_{BC}(\cdot)\chi^+(\cdot)$ where $\rho_{BC}$ is the correlation function and $\chi^+$ the characteristic function of $[0,\infty)$.  This Fourier transform, {\it i.e.},
$$	\hat\rho_{BC}^+(\omega)
	=\int_0^\infty e^{i\omega t}dt\int\rho(dx)\,B(x)C(f^{-t}x)      $$
is the boundary value on the real axis of a function of $\omega$ holomorphic for ${\rm Im}\omega>0$.  Furthermore, this function continues meromorphically to $\{\omega:{\rm Im}\omega>-\delta^*\}$ for some $\delta^*>0$, and is regular at $\omega=0$ (see [22], [26]).  Our ambition is to prove that its value at $0$ is, up to the factor $-\nu(\psi)$, equal to $(7)$.  To do this we follow the calculation in [26] Section 4 which expresses the Fourier transform as a series converging in the sense of distributions.  Up to an additive term holomorphic in $\omega$ near $\omega=0$, this is 
$$	={1\over\nu(\psi)}\tilde\nu[\tilde B_\omega\sum_{n=0}^\infty
(S^{-1}{\cal L}_{\Phi-i\omega\Psi}S)^n\tilde C_{-\omega}]\eqno{(8)}      $$
In this formula $\tilde\nu$ is the image of $\nu$ by the projection $\Sigma\to\Sigma_-$ where $\Sigma_-$ is the semi-infinite subshift defined by $\Sigma_-=\{(\xi_j)_{j\le0}:\tau_{\xi_{j-1}\xi_j}=1\}$, and $\tilde B_\omega$, $\tilde C_{-\omega}$ are H\"older continuous functions on $\Sigma_-$ depending holomorphically on $\omega$.  The {\it interactions} $\Phi$ and $\Psi$ are related to $\phi^{(u)}$ and $\psi$ and the {\it transfer operator} ${\cal L}_{\Phi-i\omega\Psi}$ acting on H\"older continuous functions on $\Sigma_-$ depends holomorphically (in fact linearly) on $\omega$.  For small $|\omega|$ the operator ${\cal L}_{\Phi-i\omega\Psi}$ is quasicompact: it has a simple eigenvalue $\lambda(\omega)$ with $\lambda(0)=1$, $\lambda'(0)\ne0$, and the rest of the spectrum is contained in a disc of radius $<1$.  The eigenfunction $S$ of ${\cal L}_{\Phi}$ to the eigenvalue $1$ is $>0$, and we have denoted by $S$ or $S^{-1}$ the multiplication or division by that function.  The derivation of the above formula is presented in [26] with slightly different notation, and one can also see that $\tilde\nu(\tilde B_0)=\tilde\nu(\tilde C_0)=0$.  We can, in the expression $(8)$, evaluate the part corresponding to the eigenvalue $\lambda(\omega)$ of ${\cal L}_{\Phi-i\omega\Psi}$.  This part is of the form $(1-\lambda(\omega))^{-1}$ times two factors, one corresponding to $\tilde B_\omega$ and the other to $\tilde C_{-\omega}$.  Both of these factors vanish at $\omega=0$ as can be seen from [26].  Since $(1-\lambda(\omega))^{-1}$ has a simple pole at $\omega=0$, the above product vanishes there.  The Fourier transform of $\rho_{BC}(\cdot)$ is thus a distribution in $\omega$ which reduces to an analytic function of $\omega$ for $|\omega|$ small, and this analytic function is given by a convergent series corresponding to the part of the spectrum of ${\cal L}_{\Phi-i\omega\Psi}$ strictly inside the unit circle.  One can thus take $\omega=0$ and obtain a convergent expression for the Fourier transform of $\rho_{BC}^+(\cdot)$ at $\omega=0$.  Manipulations as described in [26] then show that the Fourier transform of $\rho_{BC}^+(\cdot)$ at $\omega=0$ is, up to a factor $-\nu(\psi)$, equal to (7).  From Proposition 19, (6) and (7) we obtain thus Theorem B since $D_xB=D_xA$, and $\rho_{BC}=\rho_{AC}$.  
\medskip
	Note now that 
$$	\rho_{AC}(t)=\rho((A\circ f^t).C)
	=\int\rho(dx)\,A(f^tx)({\rm div}_v^{cu}(X^c+X^u))(x)      $$
$$	=-\int\rho(dx)\,(X^c(x)+X^u(x))\cdot\nabla_x(A\circ f^t)      $$
where we have used the fact that $v$ is the conditional measure of $\rho$ on center-unstable manifolds, and performed an integration by parts.  Theorem A follows then readily from Theorem B.  
\medskip
	{\bf Acknowledgments.}
\medskip
	Over the long period of preparation of this article, I have benefited from discussions with a number of colleagues, in particular V. Baladi, D. Dolgopyat, C. Liverani, and L.-S. Young.
\vfill\eject\noindent
APPENDIX
\medskip
	{\bf 1.  Calculation of $Z'$.}
\medskip
	We have
$$	f_{a*}^tx=f^tx+a{\cal R}_x^t(X^s+X^u)      $$
where we have defined, for a vector field $Y$, 
$$ {\cal R}_x^tY=\int_0^td\theta\,(T_{f^\theta x}f^{t-\theta})Y(f^\theta x) $$
Therefore 
$$	Z'=\int_0^{\psi(\xi)}dt\,(D_{f^t\pi\xi}B)
	[(T_{\pi\xi}f^t)(U^s+U^u)+{\cal R}_{\pi\xi}^t(X^s+X^u)]      $$
Notice also that 
$$	(T_{\pi\sigma^{-1}\xi}f^{\psi(\sigma^{-1}\xi)})U^{s,u}(\sigma^{-1}\xi)
+{\cal R}_{\pi\sigma^{-1}\xi}^{\psi(\sigma^{-1}\xi)}X^{s,u}=U^{s,u}(\xi)    $$
Defining
$$   ({\cal R}Y)(\xi)={\cal R}_{\pi\sigma^{-1}\xi}^{\psi(\sigma^{-1}\xi)}Y   $$
$$	({\cal T}V)(\xi)
=T_{\pi\sigma^{-1}\xi}f^{\psi(\sigma^{-1}\xi)}V(\sigma^{-1}\xi)\qquad,\qquad
	({\cal T}_-V)(\xi)=T_{\pi\sigma\xi}f^{-\psi(\xi)}V(\sigma\xi)      $$
we find 
$$   U^s=(1-{\cal T})^{-1}{\cal R}X^s=\sum_0^\infty{\cal T}^n{\cal R}X^s   $$
$$	U^u=-{\cal T}_-(1-{\cal T}_-)^{-1}{\cal R}X^u
	=-\sum_1^\infty{\cal T}_-^n{\cal R}X^u      $$
where the series on the right-hand side converge exponentially, and 
$$	(T_{\pi\xi}f^t){\cal T}^n{\cal R}X^s=(T_{\pi\sigma^{-n}\xi}
	f^{\psi(\sigma^{-n}\xi)+\cdots+\psi(\sigma^{-1}\xi)+t})
	{\cal R}_{\pi\sigma^{-n-1}\xi}^{\psi(\sigma^{-n-1}\xi)}X^s      $$
$$   =\int_0^{\psi(\sigma^{-n-1}\xi)}d\theta\,(T_{f^\theta\pi\sigma^{-n-1}\xi}
	f^{\psi(\sigma^{-n-1}\xi)+\cdots+\psi(\sigma^{-1}\xi)+t-\theta})
	X^s(f^\theta\pi\sigma^{-n-1}\xi)      $$
$$	=\int_{-\psi(\sigma^{-n-1}\xi)-\cdots-\psi(\sigma^{-1}\xi)}
	^{-\psi(\sigma^{-n}\xi))-\cdots-\psi(\sigma^{-1}\xi)}d\theta'
	(T_{f^{\theta'}\pi\xi}f^{t-\theta'})X^s(f^{\theta'}\pi\xi)      $$
Similarly
$$	(T_{\pi\xi}f^t){\cal T}_-^n{\cal R}X^u=(T_{\pi\sigma^n\xi}
	f^{-\psi(\sigma^{n-1}\xi)-\cdots-\psi(\xi)+t})
	{\cal R}_{\pi\sigma^{n-1}\xi}^{\psi(\sigma^{n-1}\xi)}X^u      $$
$$=\int_0^{\psi(\sigma^{n-1}\xi)}d\theta\,(T_{f^\theta\pi\sigma^{n-1}\xi}
	f^{-\psi(\sigma^{n-2}\xi)-\cdots-\psi(\xi)+t-\theta})
	X^u(f^\theta\pi\sigma^{n-1}\xi)      $$
$$	=\int_{\psi(\sigma^{n-2}\xi)+\cdots+\psi(\xi)}
	^{\psi(\sigma^{n-1}\xi))+\cdots+\psi(\xi)}d\theta'
	(T_{f^{\theta'}\pi\xi}f^{t-\theta'})X^u(f^{\theta'}\pi\xi)      $$
We have thus 
$$	(T_{\pi\xi}f^t)U^s+{\cal R}_{\pi\xi}^tX^s      $$
$$	=\sum_{n=0}^\infty\int_{\Psi(-n-1,\xi)}^{\Psi(-n,\xi)}d\theta\,
	(T_{f^{\theta}\pi\xi}f^{t-\theta})X^s(f^{\theta}\pi\xi)
  +\int_0^td\theta\,(T_{f^{\theta}\pi\xi}f^{t-\theta})X^s(f^{\theta}\pi\xi)  $$
$$	(T_{\pi\xi}f^t)U^u+{\cal R}_{\pi\xi}^tX^u      $$
$$	=-\sum_{n=1}^\infty\int_{\Psi(n-1,\xi)}^{\Psi(n,\xi)}d\theta\,
	(T_{f^{\theta}\pi\xi}f^{t-\theta})X^u(f^{\theta}\pi\xi)
  +\int_0^td\theta\,(T_{f^{\theta}\pi\xi}f^{t-\theta})X^u(f^{\theta}\pi\xi)  $$
$$	=-\sum_{n=2}^\infty\int_{\Psi(n-1,\xi)}^{\Psi(n,\xi)}d\theta\,
(T_{f^{\theta}\pi\xi}f^{t-\theta})X^u(f^{\theta}\pi\xi)-\int_t^{\psi(\xi)}
d\theta\,(T_{f^{\theta}\pi\xi}f^{t-\theta})X^u(f^{\theta}\pi\xi)      $$
We can also write 
$$	(T_{\pi\xi}f^t)U^s+{\cal R}_{\pi\xi}^tX^s      $$
$$	=\sum_{k=-\infty}^{-1}\int_{\Psi(k,\xi)}^{\Psi(k+1,\xi)}d\theta\,
	(T_{f^{\theta}\pi\xi}f^{t-\theta})X^s(f^{\theta}\pi\xi)
  +\int_0^td\theta\,(T_{f^{\theta}\pi\xi}f^{t-\theta})X^s(f^{\theta}\pi\xi)  $$
$$	(T_{\pi\xi}f^t)U^u+{\cal R}_{\pi\xi}^tX^u      $$
$$	=-\sum_{k=1}^\infty\int_{\Psi(k,\xi)}^{\Psi(k+1,\xi)}d\theta\,
(T_{f^{\theta}\pi\xi}f^{t-\theta})X^u(f^{\theta}\pi\xi)-\int_t^{\psi(\xi)}
d\theta\,(T_{f^{\theta}\pi\xi}f^{t-\theta})X^u(f^{\theta}\pi\xi)      $$
These two formulas give the desired evaluation of $Z'$.
\medskip
	{\bf 2.  Calculation of $\nu(Z'')$.}
\medskip
	We have 
$$	\int\nu(d\xi)\int_0^{\psi(\xi)}dt\,
	[\eta(f^t\pi\xi)-\eta(\pi\sigma^{-n}\xi)]B(f^t\pi\xi)      $$
$$	=\int\nu(d\xi)\int_0^{\psi(\xi)}dt\,B(f^t\pi\xi)
\int_{\Psi(-n,\xi)}^td\theta\,{d\over d\theta}\eta(f^\theta\pi\xi)      $$
Using charts where ${\cal X}$ is the unit vector in the last coordinate direction, we see that 
$${d\over d\theta}\eta(f^\theta\pi\xi)=({\rm div}_v^{cu}X^c)(f^\theta\pi\xi)$$
Since $\int\nu(d\xi)\eta(\pi\sigma^{-n}\xi)\int_0^{\psi(\xi)}dt\,B(f^t\pi\xi)$ tends to $0$ for $n\to\infty$ (by exponential decay of correlations for $(\nu,\sigma)$) we have 
$$	\nu(Z'')
	=\int\nu(d\xi)\int_0^{\psi(\xi)}dt\,\eta(f^t\pi\xi)B(f^t\pi\xi)      $$
$$	=\lim_{n\to\infty}\int\nu(d\xi)\int_0^{\psi(\xi)}dt\,B(f^t\pi\xi)
\int_{\Psi(-n,\xi)}^td\theta\,({\rm div}_v^{cu}X^c)(f^\theta\pi\xi)      $$
$$	=\sum_{k=-\infty}^{-1}\int\nu(d\xi)\int_0^{\psi(\xi)}dt\,B(f^t\pi\xi)
	\int_{\Psi(k,\xi)}^{\Psi(k+1,\xi)}d\theta\,
	({\rm div}_v^{cu}X^c)(f^\theta\pi\xi)      $$
$$	+\int\nu(d\xi)\int_0^{\psi(\xi)}dt\,B(f^t\pi\xi)\int_0^td\theta\,
	({\rm div}_v^{cu}X^c)(f^\theta\pi\xi)      $$
\vfill\eject\noindent
REFERENCES.
\medskip

[1] D.V. Anosov.  ``Geodesic flows on compact Riemann manifolds of negative curvature.''  Proc. Steklov Inst. Math. {\bf 90},1-209(1967).

[2] V.I. Bakhtin.  ``Random processes generated by a hyperbolic sequence of mappings. I.  Russian Acad. Sci. Izv. Math. {\bf 44},247-279(1995).  

[3] C. Bonatti, L. Diaz, and M. Viana.  {\it Dynamics beyond uniform hyperbolicity: a global geometric and probabilistic approach.} Springer, to appear.

[4] R. Bowen.  ``Markov partitions for Axiom A diffeomorphisms.''  Amer. J. Math. {\bf 92},725-747(1970).

[5] R. Bowen.  ``Periodic orbits for hyperbolic flows.''  Amer. J. Math. {\bf 94},1-30(1972).

[6] R. Bowen.  {\it Equilibrium states and the ergodic theory of Anosov diffeomorphisms.} Lecture Notes in Math. {\bf 470}, Springer, Berlin,
1975.

[7] R. Bowen.  {\it On Axiom A diffeomorphisms.}  CBMS Regional Conference No 35, Amer. Math. Soc., Providence R.I., 1978.

[8] R. Bowen and D.Ruelle.  ``The ergodic theory of Axiom A flows.''  Invent. Math. {\bf 29},181-202(1975).  

[9] N. Chernov.  ``Markov approximations and decay of correlations for Anosov flows.''  Annals of Math. {\bf 147},269-324(1998).  

[10] G. Contreras.  ``Regularity of topological and metric entropy of hyperbolic flows.''  Math. Z. {\bf 210},97-111(1992).

[11] G. Contreras.  ``The derivatives of equilibrium states.''  Bol. Soc. Brasil. Mat. (N.S.) {\bf 26},211-228(1995).  

[12] D. Dolgopyat.  ``Decay of correlations in Anosov flows.''  Annals of Math. {\bf 147},357-390(1998)' 

[13] D. Dolgopyat.  ``Prevalence of rapid mixing in hyperbolic flows.''  Ergod. Th. and Dynam. Syst. {\bf 18},1097-1114(1998).  ``Prevalence of rapid mixing-II: topological prevalence''  Ergod. Th. and Dynam. Syst. {\bf 20},1045-1059(2000).

[14] D. Dolgopyat.  ``On differentiability of SRB states for partially hyperbolic systems.'' Inventiones Math. {\bf 155},389-449(2004).

[15] M. Field, I. Melbourne, A. T\"or\"ok.  ``Stability of mixing for hyperbolic flows.''  Preprint.

[16] G. Gallavotti and E.G.D. Cohen.  ``Dynamical ensembles in stationary states.'' J. Statist. Phys. {\bf 80},931-970(1995).

[17] A. Katok, G. Knieper, M. Pollicott, and H. Weiss.  ``Differentiability and analyticity of topological entropy for Anosov and geodesic flows.''  Invent. Math. {\bf 98},581-597(1989).

[18] F. Ledrappier and J.-M. Strelcyn.  ``A proof of the estimation from below in Pesin's entropy formula.''  Ergod. Th. and Dynam. Syst. {\bf 2},203-219(1982).

[19] F. Ledrappier and L.S.Young.  ``The metric entropy of diffeomorphisms: I. Characterization of measures satisfying Pesin's formula.  II. Relations between entropy, exponents and dimension.''  Ann. of Math. {\bf 122},509-539,540-574(1985).

[20] C. Liverani.  ``On contact Anosov flows.''  Preprint.

[21] R. de la Llave, J.M. Marco and R. Moriyon.  ``Canonical perturbation theory of Anosov systems and regularity results for the Livsic cohomology equation.''  Ann. of Math. {\bf 123},537-611(1986).

[22] M. Pollicott.  ``On the rate of mixing of Axiom A flows.''  Invent. Math. {\bf 81},423-426(1982).  

[23] M. Ratner.  ``Markov partitions for Anosov flows on 3-dimensional manifolds.''  Mat. Zam. {\bf 6},693-704(1969).

[24] D.Ruelle.  ``A measure associated with Axiom A attractors.''  Am. J. Math. {\bf 98},619-654(1976).

[25] D. Ruelle.  {\it Thermodynamic formalism.}  Addison-Wesley, Reading (Mass.),1978.

[26] D. Ruelle.  ``Resonances for Axiom A flows.''  J. Differential Geometry. {\bf 25},99-116(1987).

[27] D. Ruelle.  I. ``Differentiation of SRB states.''; II. ``Correction and complements.''  Commun. Math. Phys. {\bf 187},227-241(1997); {\bf 234},185-190(2003).

[28] D. Ruelle.  ``Smooth dynamics and new theoretical ideas in nonequilibrium statistical mechanics.''  J. Statist. Phys. {\bf 95},393-468(1999).

[29] Ya.G. Sinai.  ``Markov partitions and C-diffeomorphisms.''  Funkts. Analiz i Ego Pril. {\bf 2}, No {\bf 1},64-89(1968).  English translation, Functional Anal. Appl. {\bf 2},61-82(1968).

[30] Ya.G. Sinai.  ``Constuction of Markov partitions.''  Funkts. Analiz i Ego Pril. {\bf 2}, No {\bf 3},70-80(1968).  English translation, Functional Anal. Appl. {\bf 2},245-253(1968).

[31] Ya.G. Sinai.  ``Gibbsian measures in ergodic theory.''  Uspehi Mat. Nauk {\bf 27}, No {\bf 4},21-64(1972).  English translation, Russian Math. Surveys {\bf 27}, No {\bf 4},21-69(1972).

[32] S. Smale.  ``Differentiable dynamical systems.''  Bull. AMS {\bf 73},747-817(1967).  

[33] L.-S. Young.  ``Statistical properties of dynamical systems with some hyperbolicity.''  Ann. of Math. {\bf 147},585-650(1998).

\end